\documentclass[12pt]{amsart}
\usepackage{amssymb}

\newcommand{\Z}{{\mathbb Z}}
\newcommand{\G}{{\mathbb G}}
\newcommand{\C}{{\mathbb C}}

\newcommand{\Q}{{\mathbb Q}}
\newcommand{\bP}{{\mathbb P}}
\newcommand{\A}{{\mathbb A}}
\newcommand{\cO}{{\mathcal O}}
\newcommand{\cE}{{\mathcal E}}
\newcommand{\cF}{{\mathcal F}}

\newcommand{\BoldP}{{\mathbf P}}

\newcommand{\toto}{\operatornamewithlimits{\rightrightarrows}}

\newcommand{\eqto}{\stackrel{\lower1.5pt\hbox{$\scriptstyle\sim\,$}}\to}
\DeclareMathOperator{\Hom}{Hom}
\DeclareMathOperator{\Obj}{Obj}
\DeclareMathOperator{\Spec}{Spec}

\DeclareMathOperator{\Br}{Br}

\DeclareMathOperator{\sgn}{sgn}
\DeclareMathOperator{\Sch}{Sch}
\DeclareMathOperator{\Iso}{Iso}
\DeclareMathOperator{\sets}{sets}
\DeclareMathOperator{\spec}{Spec}
\DeclareMathOperator{\GL}{GL}
\DeclareMathOperator{\PGL}{PGL}
\DeclareMathOperator{\Out}{Out}
\DeclareMathOperator{\Aut}{Aut}
\DeclareMathOperator{\End}{End}
\DeclareMathOperator{\Gal}{Gal}
\DeclareMathOperator{\band}{Band}
\theoremstyle{plain}
\newtheorem  {thm}        {Theorem}     [section]
\newtheorem  {lemma}[thm] {Lemma}
\newtheorem  {prop} [thm] {Proposition}

\newtheorem  {cor}  [thm] {Corollary}

\newtheorem* {thm*}       {Theorem}
\newtheorem* {prop*}      {Proposition}

\theoremstyle{definition}
\newtheorem  {defn} [thm] {Definition}
\newtheorem  {exml} [thm] {Example}

\theoremstyle{remark}
\newtheorem {remark} [thm]  {Remark}

\begin{document}

\title{
  Brauer groups and quotient stacks
}
\author[D. Edidin]{Dan Edidin}
\address{
  Department of Mathematics,
  University of Missouri,
  Columbia, MO 65211
}
\email{edidin@math.missouri.edu}
\author[B. Hassett]{Brendan Hassett}
\address{
  Department of Mathematics--MS 136,
  Rice University,
  6100 S. Main St.,
  Houston, TX 77005-1892
}
\email{hassett@rice.edu}
\author[A. Kresch]{Andrew Kresch}
\address{
  Department of Mathematics,
  University of Pennsylvania,
  Philadelphia, PA 19104
}
\email{kresch@math.upenn.edu}
\author[A. Vistoli]{Angelo Vistoli}
\address{
  Dipartimento di Matematica,
  Universit\`a di Bologna,
  40126 Bologna, Italy
}
\email{vistoli@dm.unibo.it}

\thanks{Edidin received support from the NSA, NSF, and the
University of Missouri Research Board while preparing this paper.
Hassett and Kresch were  partially supported by NSF Postdoctoral 
Research Fellowships. Hassett  received additional support from 
the Institute of  Mathematical Sciences of the Chinese University 
of Hong Kong and NSF. Vistoli was partially supported by
the University of Bologna, funds for selected research topics.}

\maketitle
\begin{abstract}
A natural question is to determine which algebraic stacks are
qoutient stacks. In this paper we give some partial answers and
relate it to the old question of whether, for a scheme X, the
natural map from the Brauer goup (equivalence classes of Azumaya
algebras) to the cohomological Brauer group (the torsion
subgroup of $H^2(X,{\mathbb G}_m)$ is surjective.
\end{abstract}

\section{Introduction} \label{s.intro}

Quotients of varieties by algebraic groups arise in many situations,
for instance in the theory of moduli,
where moduli
spaces are often naturally constructed as quotients of parameter spaces by
linear algebraic groups.
The quotient of a scheme by a group need not exist as
a scheme (or even as an algebraic space), and even when a quotient
exists,
the quotient morphism may not have expected
properties. For example, if $Z$ and $G$ are smooth, then the morphism
$Z \rightarrow Z/G$ need not be smooth.

To overcome this difficulty, it is often helpful
to consider quotients as {\it stacks}, rather than as schemes or algebraic
spaces.  If $G$ is a flat group scheme acting on an algebraic space $Z$
($G$ must be separated and finitely presented over some base scheme,
with the space $Z$ and the action map defined over this base),
then a quotient $[Z/G]$ always exists as a stack,
and this stack is algebraic.
Knowing that an algebraic stack
has a presentation as a quotient $[Z/G]$ (with $G$ a linear algebraic
group, say) can
make the stack easier to study, for then
the geometry of the stack is the $G$-equivariant
geometry on the space $Z$.

A natural question is to determine which algebraic stacks are quotient stacks.
In this paper we give some partial answers to this question and relate
it to the old question
of whether, for a scheme $X$, the natural map from the Brauer group
(classes of Azumaya algebras modulo an equivalence relation)
to the cohomological Brauer
group (the torsion subgroup of \'etale $H^2(X,\G_m)$) is surjective.

Some quick answers to this natural question are (the first two are
folklore):
(i) all orbifolds are quotient stacks (Theorem \ref{thm.dm});
(ii) all regular Deligne-Mumford stacks of dimension $\le 2$ are
quotient stacks
(Example \ref{e.dimensiontwo});
(iii) there exists a Deligne-Mumford stack, normal and of finite type
over the complex numbers (but singular and nonseparated) which is
not a quotient stack (Example \ref{e.noquotient}).

In fact, the example in (iii) is a stack
with stabilizer group $\Z/2$ at every point;
it is a gerbe over a normal (but nonseparated) scheme, of dimension 2
over the complex numbers.
Theorem \ref{t.brauer} says such a stack is a quotient stack
if and only if a certain class in the cohomological Brauer group
associated with it lies in the image of the map from the Brauer group. So,
(iii) yields an example, of independent interest, of non-surjectivity of
the Brauer map for a finite-type, normal, but nonseparated scheme
(Corollary \ref{c.brauernonsurject}). This stands in contrast with the
recent result of S.~Schr\"oer \cite{schr}, which says that
the Brauer map is surjective for any separated geometrically normal
algebraic surface.

The paper is organized as follows:
In Section \ref{s.stacks} we review the definition of algebraic stacks
and state accompanying results relative to quotient stacks.
Additional results concern finite covers of stacks by schemes.
In Section \ref{s.stacksgerbes} we review gerbes and Brauer groups
and state the result relating the Brauer map to
gerbes being quotient stacks. Finally in Section \ref{s.proofs}
we give proofs.

{\bf Acknowledgements.}
The authors thank Andrei Caldararu, Bill Graham, and
Amnon Yekutieli for helpful discussions. They are also
grateful to Laurent Moret-Bailly and the referee for
a number of corrections and suggestions.

\section{Stacks and quotient stacks} \label{s.stacks}

\subsection{Stacks}
Here we give a brief review of stacks. Some references
are \cite{D-M}, \cite{Vi} and \cite{L-MB}.

Stacks are categories fibered in groupoids satisfying descent-type
axioms; the stacks of interest to us will be algebraic and hence
admit descriptions in the form of groupoid schemes.
First, recall that a groupoid is
a small category $C$ in which all arrows are isomorphisms.
Write $R = \Hom(C)$ and $X = \Obj(C)$.  There are two maps $s,t\colon R
\rightarrow X$ sending a morphism to its source and target,
respectively; a map $e\colon X \rightarrow R$ taking an object to the
identity morphism of itself; a map $i\colon R \rightarrow R$ taking a
morphism to its inverse; and a map $m\colon R \times_{t,X,s} R\rightarrow R$
taking a pair of composable morphisms to their composition.
Write $j = (t,s)\colon R \rightarrow X \times X$.
There are obvious compatibilities between these maps.

A groupoid scheme consists of schemes $R$ and $X$ defined
over a fixed base scheme $L$, together with maps
$s,t,e,i,m$ satisfying the same compatibility conditions as above.
A groupoid scheme is called \'etale (respectively smooth, respectively
flat) if the maps $s$ and $t$ are
\'etale (resp.\ smooth, resp.\ faithfully flat
and locally of finite presentation).  
The stabilizer of a groupoid scheme is the scheme $S =
j^{-1}(\Delta_X)$ (here $\Delta_X \subset X \times X$ is the
diagonal).  This is a group scheme over $X$.

Let $L$ be a fixed ground scheme and let $F$ 
be a category together with a functor $p\colon F \rightarrow \Sch/L$.
For a fixed $L$-scheme $B$, let $F(B)$ denote the subcategory of $F$
consisting of objects mapping to $B$ and morphisms mapping to $1_B$.
Roughly,  a {\it category fibered in groupoids} (over
$L$) is pair consisting of a category $F$ and a
functor $p\colon F\rightarrow \Sch/L$, such that: \\
(i) For all $L$-schemes $B$, $F(B)$ is a groupoid. \\
(ii) 
For any morphism of $L$-schemes $f\colon B' \rightarrow B$
and any object $x \in F(B)$,
there is an object $f^*x$ in $F(B')$, unique up to canonical isomorphism,
together with a morphism $f^*x \rightarrow x$ lying over $f$.
For the precise definition see, e.g., \cite[Sec.\ 4]{D-M}.

A morphism of categories fibered in groupoids is simply 
a functor commuting with the projection functors
to $\Sch/L$. An isomorphism of categories fibered in groupoids
is a morphism which is an equivalence of categories.

Any contravariant functor $\Sch/L \rightarrow \sets$ determines
a category fibered in groupoids. We say that a category fibered in
groupoids over $L$ is represented
by a scheme (resp.\ algebraic space) if it is equivalent to
the functor of points of a scheme (resp.\ algebraic space).
 
An important construction is the fiber product.
Given morphisms $f_1\colon F_1 \rightarrow F$ and $f_2\colon F_2
\rightarrow F$, the fiber product $F_1 \times_F F_2$ is the category
fibered in groupoids defined as follows:
Objects are triples $(x_1, x_2, \psi)$ where
$x_1$ is an object of $F_1$, $x_2$ is an object of
$F_2$, and $\psi\colon f_1(x_1) \rightarrow f_2(x_2)$ is an isomorphism,
lying over an identity morphism of $\Sch/L$.
A morphism is specified by a pair of morphisms
compatible with the induced isomorphism in $F$.

\begin{defn}
A category fibered in groupoids $(F,p)$ is a {\it stack}
if it satisfies
two descent properties.\\
(1) For objects $x$, $y$ in $F(B)$ the functor
$\Iso_B(x,y)\colon \Sch/B \rightarrow \sets$ assigning
to a $B$-scheme $f\colon B' \rightarrow B$
the set of isomorphisms between $f^*x$ and $f^*y$ is a sheaf
for the \'etale topology.\\
(2) $F$ has effective descent
for \'etale morphisms.
\end{defn}

\begin{defn}
A morphism of stacks
is {\it representable} if for any morphism of an algebraic
space
$B \rightarrow F$, the fiber product $B \times_F F'$
is represented by an algebraic space. A morphism
is {\it strongly representable} if for any morphism
of a scheme $B \rightarrow F$, the fiber product
$B \times_F F'$ is represented by a scheme.
\end{defn}

Let $\BoldP$ be a property of morphism of schemes which is preserved
by base change  and is local for the smooth topology. A representable
morphism $F' \rightarrow F$ has property $\BoldP$ if for
all morphisms $B \rightarrow F$ of algebraic spaces, the induced
morphism $B \times_F F' \rightarrow B$ has property $\BoldP$.

Stein factorization holds for algebraic spaces and
implies \cite[II.6.15]{Kn}
that if $f\colon X\rightarrow Y$ is a separated quasi-finite morphism
of algebraic spaces, and if $Y$ is a scheme, then $X$ is a scheme.
Hence, a representable separated quasi-finite morphism is
always strongly representable.

\begin{defn} \label{d.algstack}
A stack $F$ is {\it algebraic}, or is an {\it Artin stack}, if \\
(1) There exists a representable smooth surjective morphism 
$X \rightarrow F$ from a scheme. \\
(2) The diagonal morphism $F \rightarrow F \times_L F$
is representable, quasi-compact, and separated.
\end{defn}

\begin{remark}
The representability of the diagonal implies
that any morphism from an algebraic space is representable.
For stacks with quasi-finite diagonal,
any morphism from a scheme is strongly representable.
\end{remark}

\begin{remark} \label{rem.dmstack}
A stack $F$ is called a {\it Deligne-Mumford stack} 
if there exists an \'etale cover of $F$ by a scheme. By \cite[8.1]{L-MB},
this happens if and only if the diagonal $F \to F \times_L F$ is
unramified\footnote{As explained in \cite[4.2]{L-MB}, `unramified' should
be understood to mean `locally of finite type and formally unramified'.}.
A Deligne-Mumford
stack has, in particular, quasi-finite diagonal.
The geometric fibers of the diagonal are group schemes,
so if all the residue fields of $L$ 
have characteristic 0 then conversely, any
algebraic stack with quasi-finite diagonal
is a Deligne-Mumford stack.
\end{remark}

Finally, we describe very briefly groupoid presentations (or atlases) of
algebraic stacks: see \cite{L-MB} for a full treatment. By definition,
any algebraic stack
$F$ admits a smooth surjective map from a scheme $X$; $X \rightarrow F$ is
called a {\em smooth atlas}.
In this case, the fiber product $R = X \times_{F} X$
is an algebraic space. However, for stacks with quasi-finite
diagonal, the diagonal is strongly representable, so
$R$ is in fact a scheme. The smooth
groupoid scheme $R \toto X$ is called a {\em presentation} for $F$.
Conversely, any smooth groupoid scheme $R \toto X$ with
separated, finite-type relative diagonal $R\to X\times X$
determines an algebraic stack $[R \toto X]$.
A theorem of  Artin (cf.\ \cite[10.1]{L-MB}) says that any
faithfully flat groupoid scheme $R \toto X$ (with separated, finite-type
relative diagonal) determines an algebraic stack. In this
case the groupoid scheme
$R \toto X$ is called a {\em faithfully flat presentation} for $F$.
By Remark \ref{rem.dmstack}, an algebraic stack
is a Deligne-Mumford stack if and only if
it has an {\em \'etale presentation}.
If $F$ is an algebraic stack with quasi-finite diagonal, then it is relatively
straightforward (\cite[Lemma 3.3]{KM}) to show that
$F$ has a quasi-finite faithfully flat atlas of schemes.

If the group scheme $G$ acts on the algebraic space $Z$
(we assume $G$ flat, separated, and of finite presentation over the
ground scheme $L$; the space $Z$ should be an $L$-space and the action
map $Z\times_L G\to Z$ an $L$-morphism), then the action determines a
groupoid $Z\times_L G\toto Z$.
This will be a flat atlas for the stack whose fiber over any
$L$-scheme $T$ is the category of principal $G$-bundles $E\to T$
together with $G$-equivariant morphisms $E\to Z$.
This is an algebraic stack, denoted $[Z/G]$.

As noted above, any algebraic space is an algebraic stack; the following
result says when the converse holds.

\begin{prop}(\cite[2.4.1.1 and 10.1]{L-MB}) \label{algspace}
Let $F$ be an algebraic stack, and let
$s,t\colon R\toto X$ be a faithfully flat presentation for $F$.
Then $F$ is an algebraic
space if and only if the map $R \rightarrow X \times X$ is a
monomorphism.
If we set $S = (t \times s)^{-1}(\Delta_X)$ then this is equivalent
to 
$S \rightarrow X$ being an isomorphism
by either $s$ or $t$.
\end{prop}

We call attention to the map $S\to X$ of Proposition \ref{algspace}.
The fiber product of the diagonal $F\to F\times F$ with itself
is an algebraic stack
$I_F:=F\times_{F\times F}F$.
The projection (to either factor) $I_F\to F$
is the {\em stabilizer map}, and is represented by the
stabilizer $S\to X$ of the groupoid space $R\toto X$,
for any atlas $X$.

\subsection{Results on stacks} \label{s.stackresults}

The first theorem states that stacks with quasi-finite diagonal
are finitely para\-met\-rized; i.e., admit finite covers by schemes.
This is the strongest possible result since any finitely parametrized
stack must have quasi-finite diagonal.
This result extends results of
Vistoli \cite{Vi} and Laumon and Moret-Bailly \cite{L-MB} for Deligne-Mumford
stacks. The first
result of this form of which the authors are aware is due to 
Seshadri \cite[Theorem 6.1]{Se} in the context
of group actions on varieties. In fact, the use of Lemma \ref{l.shtrick}
was inspired by reading his paper.

\begin{thm} \label{thm.fp}
Let $F$ be an algebraic stack
of finite type over a Noetherian ground scheme $L$.
Then the diagonal $\delta\colon F \rightarrow F \times_L F$ is quasi-finite
if and only
if there exists a finite
surjective morphism $X \rightarrow F$ from a (not necessarily separated) 
scheme $X$.
\end{thm}

\begin{remark}
Existence of finite scheme covers is an important ingredient in
intersection theory on Deligne-Mumford stacks.
It is used, for instance, to
define proper pushforward for nonrepresentable morphisms
(of cycles modulo rational equivalence with $\Q$ coefficients).
General intersection-theoretic machinery has recently been developed for
Artin stacks whose geometric stabilizers are affine groups \cite{Kr}.
All of intersection theory on Deligne-Mumford stacks, as in
\cite{Gi} and \cite{Vi},
generalizes to
Artin stacks with quasi-finite diagonal, where
Theorem \ref{thm.fp} is used to provide nonrepresentable proper pushforwards.
\end{remark}

\begin{defn}
\label{defnquotient}
Let $F$ be a stack, of finite type over a Noetherian base scheme $L$.
We say $F$ is a {\em quotient stack} if $F$ is isomorphic to a
stack of the form $[Z/G]$ where $Z$ is an algebraic space,
of finite type over $L$, and $G$ is
a subgroup scheme of the general linear group scheme
$GL_{n,L}$ for some $n$, with $G$ flat over $L$.
\end{defn}

\begin{remark}
Every affine group scheme of finite type over a field is a subgroup
scheme of $GL_n$, so the condition on $G$ in Definition \ref{defnquotient}
is the natural notion of linear algebraic group over a general
Noetherian base.
\end{remark}

\begin{remark}
The quotient $Z':=Z\times_L GL_{n,L}/G$ (where $G$ acts on $Z$
and acts by translation on $GL_{n,L}$) exists
as an algebraic space, and
$[Z'/GL_n]\simeq[Z/G]$.
So every quotient stack is a quotient by $GL_n$ for some $n$.
\end{remark}

We state two foundational results, followed by two results giving sufficient
conditions for a stack to be a quotient stack.
Recall that $f\colon E\to F$ is a projective morphism if and only if
$f$ factors, up to $2$-isomorphism, as a closed immersion followed by
projection $E\to {\mathbf P}(\cE)\to F$, where $\cE$ is a
finite-type quasi-coherent sheaf on $F$ and ${\mathbf P}(\cE)$ denotes
its projectivization.

\begin{lemma} \label{l.vbqs}
Let $F$ be an algebraic stack of finite type over a Noetherian scheme.
The following are equivalent.
\begin{itemize}
\item[(i)] $F$ is a quotient stack.
\item[(ii)] There exists a vector bundle $V\to F$ such that at every geometric
point, the stabilizer action on the fiber is faithful.
\item[(iii)] There exists a vector bundle $V\to F$ and a locally closed substack
$V^0\subset V$ such that $V^0$ is representable and $V^0$ surjects onto $F$.
\end{itemize}
\end{lemma}

\begin{lemma} \label{l.flat}
Let $\pi\colon E\to F$ be a flat projective map of stacks
(of finite type over a Noetherian base scheme) which is surjective.
If $E$ is a quotient stack, then so is $F$.
\end{lemma}

\begin{thm} \label{thm.coq}
Let $F$ be an algebraic stack of finite type over a Noetherian scheme, and let
$f\colon X \rightarrow F$ be a finite cover by 
a scheme or algebraic space. 
If the coherent sheaf
$f_*{\mathcal O}_X$ is the quotient of a locally free coherent sheaf
then $F \simeq [Z/\GL_n]$ where $Z$ is an algebraic space.
In particular, if every coherent sheaf on $F$ is the quotient of
a locally free coherent sheaf, then $F$ is a quotient stack.
\end{thm}

\begin{remark}
If the ground scheme $L$ is normal and separated, and has the property
that every coherent sheaf on $L$ is the quotient of
a locally free sheaf (e.g., if $L$ is affine, or regular) and
if $F = [Z/\GL_n]$ where $Z$ is a scheme equivariantly
embedded in a regular Noetherian separated scheme, then the equivariant
resolution theorem of \cite{Th} implies 
that every coherent sheaf on $F$
is the quotient of a locally free coherent sheaf on $F$.
\end{remark}

\begin{cor} \label{c.ffc}
Let $F$ be an algebraic stack of finite type over a Noetherian scheme.
If $F$ has a finite flat cover by an algebraic space
then $F$ is a quotient stack.
In particular,
if $F$ is regular and has a finite cover by a Cohen-Macaulay
algebraic space then $F$ is a quotient stack.
\end{cor}

\begin{exml}
\label{e.dimensiontwo}
Assume the base scheme is a field
(or more generally, any universally Japanese scheme,
for instance $\spec \Z$).
Any regular stack of dimension $\le 2$ with quasi-finite diagonal
is a quotient stack.
\end{exml}

We emphasize the second statement of Corollary \ref{c.ffc}
because the Deligne-Mum\-ford stacks considered by Mumford
in \cite{Mumenum} satisfy (ii). In particular they are quotient
stacks, and the intersection product he constructs
is a special case of the intersection product of
\cite{E-G}.

Finally, if $F$ is Deligne-Mumford then we have the following result
which we obtained based on conversations with Bill Graham.
In the characteristic zero setting, this result is familiar from the
study of orbifolds.
\begin{thm} \label{thm.dm}
If $F$ is a smooth Deligne-Mumford stack
of finite type over the Noetherian base scheme
such that the automorphism group of 
a general geometric point of $F$ is trivial, then
$F$ is a quotient stack.
\end{thm}
Thus any stack which admits a representable morphism to a
smooth Deligne-Mumford stack with trivial generic stabilizers, also, is a
quotient stack.

\begin{cor} \label{c.trivcenter}
Let $F$ be a smooth Deligne-Mumford stack of finite type
over a Noetherian base scheme.
Assume $F$ has finite stabilizer,
and suppose the automorphism group of a general geometric point of $F$
has trivial center.
Then $F$ is a quotient stack.
\end{cor}

Recall that an algebraic space $Q$ has
quotient singularities if locally in the \'etale topology
$Q$ is isomorphic to quotients $U/H$, where $H$ is a finite group
and $U$ is smooth.
By \cite[Proposition 2.8]{Vi},
any separated scheme of finite type over a field of characteristic zero
with quotient singularities is a moduli space
for a smooth stack $F$ which has generically trivial stabilizer,
so we have the following consequence.
\begin{cor} \label{c.quot}
Any separated scheme of finite type over a field of
characteristic 0 which has at worst  quotient
singularities is a quotient $Q = Z/G$ where
$Z$ is a smooth algebraic space and $G$ is a linear algebraic group.
\end{cor}

Lastly, as promised, not every Deligne-Mumford stack is a
quotient stack.

\begin{exml} \label{e.noquotient}
Let $Y$ be the scheme $\spec \C[x,y,z]/(xy-z^2)$, whose nonsingular
locus is $Y^{\rm reg}=Y\smallsetminus\{0\}$.
There is a unique (up to 2-isomorphism) nontrivial involution of
$Y^{\rm reg}\times B(\Z/2)$ which commutes with the projection
map to $Y^{\rm reg}$.
Let $F$ be the stack gotten by glueing two copies of
$Y\times B(\Z/2)$ via this involution.
Then $F$ is not isomorphic to $[Z/G]$ for any algebraic space $Z$ and
algebraic group $G$.
\end{exml}

\section{Gerbes and Brauer groups} \label{s.stacksgerbes}
In this section we give a brief review of gerbes and
Brauer groups and 
state our accompanying results.
References for gerbes are \cite{Mi} and \cite{L-MB}.
For Brauer groups, see \cite{Gr} and \cite{Mi}.

\subsection{Gerbes}
In what follows we fix a base scheme $X$, assumed Noetherian, and we
take $G$ to be a group scheme, flat, separated, and of finite type over $X$.
The gerbes that arise in the theorem that relates gerbes to Brauer groups
(Theorem \ref{t.brauer}) have $G$ equal to the
algebraic torus $\G_m$ or a group of roots of unity $\mu_n$.
We only discuss gerbes that are modeled on some group scheme $G$ over
the base.

\begin{defn}
A {\it $G$-gerbe} over $X$ is a morphism $F\to X$, with
$F$ an algebraic stack,
such that there exists a faithfully flat map, locally of finite presentation,
$X'\to X$, such that $F\times_XX'\simeq BG\times_XX'$.
\end{defn}

We say the $G$-gerbe $F\to X$ is {\em trivial} if $F\simeq BG$.
Note that a gerbe $F\to X$ which admits a section $x\in F(X)$ satisfies
$F\simeq B(\Aut_F(x))$ where $\Aut_F(x)$ is group scheme
(or group space)
$\Iso_F(x,x)$
(such a gerbe is called {\em neutral}).
Nontrivial gerbes are easy to construct, much the way one constructs
nontrivial vector bundles, or torsors.
For instance, one can glue two copies of $\A^1\times B(\Z/2)$ along a
nontrivial involution of $(\A^1\smallsetminus\{0\})\times B(\Z/2)$
to obtain a nontrivial $(\Z/2)$-gerbe over $\bP^1$.

\begin{defn}
Let $G$ and $H$ be two group schemes over $X$.
The sheaf of {\it band isomorphisms}, denoted
$\band(G, H)$, is the sheafification of the quotient of the sheaf of group
isomorphisms $\Iso(G, H)$ by the conjugation action of $H$.
When $G=H$, this is the sheaf of {\it outer automorphisms} of $G$,
which is denoted $\Out(G)$.
\end{defn}

%$\band(G, H)$ is the set of isomorphisms of $G$ with $H$ in the category of
%bands over $X$ (we use {\it band\/} for the French word {\it lien}.) See
%\cite[??? ]{Gir}.
%
%Obviously if $H$ is abelian then $\band(G, H) = \Iso(G, H)$.

\begin{defn}
Given a $G$-gerbe $F\to X$,
the {\it associated torsor of outer automorphisms}
is the sheaf $P$ over $X$ defined as follows.
Let $T$ be an $X$-scheme.
If there exists an object $t\in F(T)$, then we define
$P(T)$ to be
$\band(\Aut_F(t), G\times_XT)$. One checks that if $\tilde t$ denotes
another object in
$F(T)$, then there is a canonical element $\band(\Aut_F(t),\Aut_F(\tilde
t\,))$ obtained by chosing local isomorphisms of $t$ with $\tilde t$; this
canonically identifies $\band(\Aut_F(t), G\times_XT)$ with
$\band(\Aut_F(\tilde t\,), G\times_XT)$.
In general, $P(T)$ is defined as the difference kernel
$P(T')\toto P(T'\times_TT')$ with respect to any
flat cover $T'\to T$ such that $t'\in F(T')$ exists.
%
%In general, one defines $P(T)$ by choosing a flat cover $T'\to T$
%such that $t'\in F(T')$ exists and setting $P(T)$ to be the
%difference kernel of $P(T')\toto P(T'\times_TT')$.
Elements of $P(T)$ pull back in the obvious fashion.
\end{defn}

There is an obvious action of $\Out(G)$ on $P$, making it into a torsor.
This torsor is
classified by some $\beta\in H^1(X,\Out(G))$.
This is the {\it first obstruction to triviality} of $F$. 

For the remainder of this section,
we assume that $\Out(G)$ is a finite flat group
scheme over $X$.
This is the case when (i) $G$ is finite group
(viewed as a group scheme over $\Spec\Z$ and hence over any base);
(ii) $G=\mu_n$ for any positive integer $n$;
(iii) $G=\G_m$.
Now there are two ways to remove the first obstruction to triviality
for a gerbe.
First, one can hope that $\beta$ is in the image of
$H^1(X,\Aut(G))\to H^1(X,\Out(G))$, and then use the $\Aut(G)$-cocycle
to substitute, in place of $G$, a new group scheme $G'$,
locally isomorphic to $G$.
For instance, if the symmetric group $S_3$ acts on
$\A^1\smallsetminus\{0\}$
by $\sigma\cdot z = \sgn(\sigma)z$, then
$F:=[\A^1\smallsetminus\{0\}/S_3]$ is nontrivial as a $(\Z/3)$-gerbe over
$X:=\A^1\smallsetminus\{0\}$.
Its first obstruction class is the nontrivial element of
$H^1(X,\Aut(\Z/3))=H^1(X,\Out(\Z/3))$.
Twisting, we obtain a group scheme $G'$ over $X$, and we
find in this example that $F\simeq BG'$.

The second method, which doesn't require hoping, is to pull back
to the total space of the $\Out(G)$-torsor.
So,
{\it the first obstruction to triviality vanishes upon finite flat pullback}.

Assume our $G$-gerbe has trivial first obstruction, and let a
trivialization of $P$ be fixed. If the center of $G$ is trivial,
then one can use
the stack axioms to glue local sections of $F\to X$ to get a section defined
over $X$; the cocycle condition will automatically be satisfied.
In general, the obstruction is a 2-cocycle with values in the center $Z$
of $G$.
The class $\alpha\in H^2(X,Z)$ is the {\it second obstruction
to triviality} of $F$
(this depends on the choice of trivialization of $P$; a different
choice will differ by a
global section $\gamma$ of $\Out(G)$, and
the class in $H^2(X,Z)$ resulting from
the new section is the result of $\gamma$ applied to $\alpha$ by the
obvious action of $\Out(G)$ on $Z$).

\begin{remark}
A gerbe $F$ is said to be
{\em banded} (it is becoming standard to translate as {\em band}
the French verb {\em lier})
by $G$ if the gerbe is endowed
with a global section of the associated torsor of outer automorphisms.
When $G$ is abelian, to say that $F$ is banded by $G$ is equivalent to
saying that for every $X$-scheme $U$ and object $u\in F(U)$,
there is chosen an isomorphism $G(U)\eqto \Aut_{F}(u)$,
compatible with pullbacks.
\end{remark}

If $G$ is abelian,
then by cohomological machinery, the set of isomorphism classes of
gerbes on $X$ banded by $G$ is in bijection with $H^2(X,G)$
(\cite[\S{}IV.2]{Mi}).
For $G$ finite and flat with flat center $Z$ let
$F\to X$ be a gerbe banded by $G$,
with second obstruction $\alpha\in H^2(X,Z)$,
then the $Z$-gerbe $E\to X$ associated with $\alpha$
admits a finite flat representable morphism to $F$.
So, such a $G$-gerbe is covered by a gerbe banded by $G$,
which in turn is covered by a gerbe banded by the center of $G$.

\begin{prop} \label{p.gerbes}
Let $G$ be a finite flat group scheme over $X$.
Assume that the center, $Z$, and the sheaf of outer automorphisms
$\Out(G)$ are finite and flat as well.
Let $F\to X$ be a $G$-gerbe.
Then there exists an $\Out(G)$-torsor $Y\to X$ and
a gerbe $E\to Y$ banded by $Z$, such that $E$ admits a finite flat
representable surjective morphism to $F$.
\end{prop}

\subsection{Brauer groups}
Let $X$ be a Noetherian scheme.
The {\it Brauer group} $Br(X)$ is the group of Azumaya algebras
(sheaves of algebras, \'etale-locally isomorhic to endomorphism algebras
of vector bundles), modulo the equivalence relation
$\cE\sim\cE'$ if $\cE\oplus\End(V)\simeq \cE'\oplus\End(V')$ for
some pair of vector bundles $V$ and $V'$ on $X$.
By the Skolem-Noether theorem,
the rank $n^2$ Azumaya algebras on $X$ are classified by
$H^1(X,\PGL_n)$.
The exact sequence
$$1\to \G_m\to \GL_n\to \PGL_n\to 1,$$
identifies the obstruction to a rank $n^2$ Azumaya algebra being
the endomorphism
algebra of a vector bundle as an element~--~in fact,
an $n$-torsion element~--~of the \'etale cohomology group $H^2(X,\G_m)$.
There is thus determined a homomorphism
$$Br(X)\to H^2(X,\G_m).$$
It is a fact that this homomorphism is always injective \cite[IV Th.\ 2.5]{Mi}.

The {\it cohomological Brauer group}, denoted $\Br'(X)$, is defined to be
the torsion subgroup of $H^2(X,\G_m)$.
When $X$ is regular, the map $\Br'(X)\to \Br'(k(X))$ is injective,
where $k(X)$ denotes the generic point of $X$.
For a field,  $\Br$, $\Br'$ and the full second cohomology group
agree.
It is only in the presence of singularities that the cohomological
Brauer group may differ from the full cohomology group $H^2(X,\G_m)$.

The {\it Brauer map} is the injective group homomorphism
$$\Br(X)\to \Br'(X).$$
A major question in the study of Brauer groups is:
for which schemes $X$ is the Brauer map an isomorphism?
The article \cite{Ho} identifies some classes of schemes
for which this is known.
The Brauer map
is known to be an isomorphism for abelian varieties, low-dimensional
varieties (general varieties of dimension $1$ and regular varieties of
dimension 2), affine varieties, and separated unions of two affine
varieties.
Recently the Brauer map has been shown to be an isomorphism for
separated geometrically normal algebraic surfaces \cite{schr}.
Also known in general is that if $\alpha\in \Br'(X)$ is trivialized by a
finite flat cover, then $\alpha$ lies in the image of the Brauer map.

\begin{thm} \label{t.brauer}
Let $X$ be a Noetherian scheme.
Let $\beta$ be an element of $H^2(X,\G_m)$.
The following are equivalent.
\begin{itemize}
\item[(i)] $\beta$ lies in the image of the Brauer map.
\item[(ii)] There exists a flat projective morphism of schemes
$\pi\colon Y\to X$,
surjective, such that $\pi^*\beta=0$ in $H^2(Y,\G_m)$.
\item[(iii)] The $\G_m$-gerbe with classifying element $\beta$
is a quotient stack.
\end{itemize}
Furthermore, if $n\beta=0$ and $\alpha\in H^2(X,\mu_n)$ is a pre-image
of $\beta$ under the map of cohomology coming from the Kummer sequence,
then conditions {\em (i)}, {\em (ii)}, and {\em (iii)} are equivalent to
\begin{itemize}
\item[(iv)] The $\mu_n$-gerbe with classifying element $\alpha$
is a quotient stack.
\end{itemize}
\end{thm}

\begin{remark}
Here we are writing `$G$-gerbe with classifying element $\alpha$'
(for $G=\G_m$ or $\mu_n$) to refer to a gerbe, banded by $G$, whose
second obstruction to triviality is $\alpha\in H^2(X,G)$.
Such a gerbe is defined uniquely up to isomorphism, hence the
abusive terminology `the $G$-gerbe$\ldots$'.
\end{remark}

\begin{remark}
In characteristic $p>0$ (or in mixed characteristic) the cohomology groups
above are flat cohomology groups. By \cite[III.11]{Gr}, sheaf cohomology
with values in $\G_m$, or in $\mu_n$ when $n$ in invertible, is the same in
the \'etale and flat topologies.
\end{remark}

\begin{remark}
Statements (i) and (ii) do not involve stacks, so the implications
(i) $\Leftrightarrow$ (ii) have independent interest.
One direction, (i) $\Rightarrow$ (ii), is well-known: if $\beta$ is in
the image of the boundary homomorphism $H^1(X,PGL_{m+1})\to H^2(X,\G_m)$
then pullback to the associated $\bP^m$-bundle trivializes $\beta$;
the $\bP^m$-bundle is the famous Brauer-Severi scheme.
The other direction, (ii) $\Rightarrow$ (i), seems to have been
known only as folklore, until recently.
The result now appears in the Ph.D. thesis of A.~Caldararu
\cite[Prop.\ 3.3.4]{Ca}.
\end{remark}

\begin{remark}
The question of whether a general Deligne-Mumford stack is
a quotient stack is hard
(even with strong hypothesis such as smooth and proper over a field).
But for gerbes over schemes over a field of characteristic zero,
Lemma \ref{l.flat} can be used, in conjunction with
Proposition \ref{p.gerbes} and Theorem \ref{t.brauer}, to
reduce the question to the case of 
$\mu_n$-gerbes.
Indeed, by Proposition \ref{p.gerbes}, any $G$-gerbe has a finite flat
representable cover
by an abelian group gerbe, which in turn admits a closed immersion to
a product of cyclic group gerbes.
\end{remark}

Example \ref{e.noquotient} then tells us:

\begin{cor} \label{c.brauernonsurject}
Let $X$ be the union of two copies of
$\Spec \C[x,y,z]/(xy-z^2)$, glued along the nonsingular locus.
Then the Brauer map $\Br(X)\to\Br'(X)$ is not surjective.
So the nonseparated union of two affine schemes need not have
surjective Brauer map.
\footnote{R. Hoobler has pointed out that it is possible to verify directly,
using cohomological methods, that the scheme $X$ in this statement satisfies
$\Br(X)=0$ and $\Br'(X)=\Z/2$.}
\end{cor}

Going the other way, Theorem \ref{t.brauer} provides an example of a
stack with affine (but not quasi-finite) 
diagonal of finite type over a field
which is not a quotient stack.
(Note that 
the stack in Example \ref{e.noquotient} only has quasi-affine diagonal.)

\begin{exml}
Let $X$ be a normal separated surface over a field (if one wishes, $\C$)
such that $H^2(X,\G_m)$ contains a non-torsion element $\beta$
\cite[II.1.11.b]{Gr}.
Then the $\G_m$-gerbe $F$ classified by $\beta$ has affine diagonal
and is not a quotient stack.
\end{exml}

\section{Proofs of results} \label{s.proofs}

\subsection{Finite parametrization of stacks}
Here we prove Theorem \ref{thm.fp}, which states that
that every stack with quasi-finite diagonal has a finite cover by a scheme.
We begin with an easy, but very useful lemma.

\begin{lemma} \label{l.shtrick}
Suppose that $p_1\colon F_1 \rightarrow F$ and $p_2\colon F_2 \rightarrow F$
are representable (respectively strongly representable) morphisms.
Assume that $F$ is covered by open substacks $U_1$, $U_2$
such that the fiber products $U_1 \times_F F_2$ and $F_1 \times_F U_2$ are 
representable by algebraic spaces (resp.\ schemes).
Then the 
fiber product is $F_1 \times_F F_2$
is also represented by an algebraic space (resp.\ scheme).
\end{lemma}

\begin{proof}
The inverse images of $U_1 \times_F F_2$ and $F_1 \times_F U_2$ in $F_1
\times_F F_2$ are represented by algebraic spaces (resp.\ schemes),
because $p_1$ and $p_2$ are representable (resp.\ strongly representable).
But these inverse images are open substacks which cover $F_1 \times_F F_2$.
\end{proof}

\begin{proof}[Proof of Theorem \ref{thm.fp}]
Since $F$ is finitely presented
over the ground scheme, we may assume that $F$ is obtained
by base change from a stack of finite type over $\spec \Z$. Hence
to obtain a cover we may assume that $F$ is of finite type over
$\spec \Z$.
Also, since the morphism $F_{red} \rightarrow F$ is finite and surjective 
we can assume $F$ is reduced.
By working with each irreducible component separately we can
assume $F$ is integral. Finally by normalizing we can assume 
that $F$ is normal.

Suppose that $F$ has an open cover $F^1, \ldots, F^k$ such
that $F^i$ has a finite cover by a scheme $Z^i$.
The composite morphism $Z^i \rightarrow F^i \hookrightarrow F$
is quasi-finite. Thus, by 
Zariski's Main Theorem \cite[Theorem 16.5]{L-MB}
the morphism $Z^i \rightarrow F$
factors as an open immersion followed by a finite representable 
map $Z^i \hookrightarrow
Z_i \rightarrow F$. Since $F$ is assumed to be irreducible,
the finite representable morphism 
$Z_i \rightarrow F$ has dense image so it must be surjective.
Set $Z = Z_1 \times_F Z_2 \ldots \times_F Z_k$.
The induced map $Z \rightarrow F$ is finite, representable
and has dense image, so it is surjective. Since
any finite representable morphism is strongly representable, we
can, by applying the Lemma,
conclude that $Z$ is a scheme.

Thus, to prove the theorem it suffices to prove that
$F$ has a cover by open substacks which admit finite covers
by schemes. By \cite[Lemma 3.3.1]{KM}, $F$ has a quasi-finite flat
cover by a scheme $V$. Let $V_i$ be an irreducible component
of $V$. Once again applying Zariski's Main Theorem, 
the quasi-finite morphism $V_i \rightarrow F$ factors
as $V_i \hookrightarrow F' \rightarrow F$, where the 
first map is an open immersion and the second map is
finite (and by density surjective). Replacing $F$ by $F'$
we may therefore assume that $F$ is generically a scheme.
In particular, we can assume that $F$ has a generic point $\spec K$.

Let $s,t\colon R\toto X$ be a smooth presentation for $F$. Since
we are working locally we can assume that $X$ is a normal
variety.
By \cite[Lemma 3.3.1]{KM}, the smooth cover can be refined
to
a quasi-finite flat cover by a scheme $V$ 
and the morphism $V \rightarrow X$ is the composition
of a closed immersion and an \'etale morphism. 
Again since we are working locally we may assume that $V$ is
irreducible. In particular we may also
assume that $V$ is normal.

Since the morphism $V \rightarrow F$ is quasi-finite, it is open.
Replacing $F$ by an open substack, we may assume that $V\rightarrow F$
is surjective.
Now we construct a finite cover of $F$ by a scheme.
The map $V \rightarrow F$ is generically
finite, so $K(V)$ is a finite extension of $K$ (recall that
$\spec K$ is the generic point of $F$). Let
$K'$ be a normal extension of $K$ containing $K(V)$.
Then $K'$ is Galois over a field $K''$ which is a purely
inseparable extension of $K$.
Let $F'$ be the normalization of $F$ in $K'$.
Let $U_1$ be the pre-image of $V$ in $F'$,
and for $\alpha\in \Gal(K'/K'')$ let $U_\alpha$ be the translate
of $U_1$ under the action of $\alpha$.
Each $U_\alpha$ is a scheme.
Since normalization commutes with smooth pullback
(\cite[Lemma 16.2.1]{L-MB}),
we may invoke \cite[Prop.\ V2.3.6]{Bourbaki}
to deduce that
$\Gal(K'/K'')$ acts 
transitively on the fibers of $F'/F$.
Hence the $U_\alpha$ cover $F'$, so $F'$ is a scheme which is a finite
cover of $F$.
\end{proof}

As a corollary of independent interest, we obtain Chow's Lemma for 
stacks with finite diagonal, extending \cite[Theorem 4.12]{D-M}
(a stack with quasi-finite diagonal is separated if and only
if the diagonal is finite).
\begin{cor} \label{c.chow}
Let $F$ be an algebraic stack of finite type over a Noetherian ground scheme.
If the diagonal of $F$ is finite,
then $F$ admits a proper, surjective, generically
finite morphism from a quasi-projective scheme.
\end{cor}

\subsection{Stacks which are quotient stacks}
In this section we give proofs of Lemmas \ref{l.vbqs} and \ref{l.flat},
and from these deduce Theorems \ref{thm.coq} and \ref{thm.dm}.

In Lemma \ref{l.vbqs}, the implication (i) $\Rightarrow$ (iii)
is well-known: if $F\simeq[X/G]$, let $G$ act linearly
on some affine space $\A^m$, freely on some open $U\subset \A^m$ such
that the structure map from $U$ to the base scheme is surjective.
Now we take $V^0\subset V$ to be $[X\times U/G]\subset[X\times \A^m/G]$
with the diagonal $G$-action.
Clearly, (iii) implies (ii).
If $V\to F$ is a vector bundle of rank $n$ such that at every geometric point,
the stabilizer action is faithful on the fiber,
then the stabilizer action on frames is free at every geometric point,
hence the associated frame bundle $P$ is
an algebraic space (Proposition \ref{algspace}),
and $F\simeq[P/GL_n]$.
This establishes (ii) implies (i), and we have proved Lemma \ref{l.vbqs}.

\smallskip

To prove Lemma \ref{l.flat},
let $E$ and $F$ be  finite-type stacks over a Noetherian ground scheme,
and let $\pi\colon E\to F$ be a flat, projective morphism.
Let $\cO(1)$ denote a relatively ample invertible sheaf on $E$,
and for a coherent sheaf $\cE$ on $E$,
we let $\cE(k)$ denote $\cE\otimes\cO(k)$.
We know that for $k$ sufficiently large,
we have $R^i\pi_*\cE(k)=0$ for $i>0$ and hence $\pi_*\cE(k)$ locally free
(these are local assertions, and for schemes this is well known).

Suppose $E$ is a quotient stack.
Then there is a locally free coherent sheaf $\cE$ on $E$, such that the
geometric stabilizer group actions on fibers are faithful.
Replacing $\cE$ by $\cE\oplus\cO_E$ if necessary, the stabilizer actions
on fibers of $\cE(k)$ for each $k$ will be faithful as well.
Choose $k$ such that $R^i\pi_*\cE(k)=0$ for $i>0$
and such that the natural map of sheaves $\pi^*(\pi_*\cE(k))\to\cE(k)$
is surjective.
We may also suppose $\cE(k)$ is very ample when restricted to
the fibers of $\pi$.
If we let $\cF:=\pi_*\cE(k)$, then $\cF$ is a locally free coherent
sheaf on $F$ such that the stabilizer group actions on fibers
are faithful.
Indeed, if $p\colon \spec\Omega\to F$ is a geometric point,
with stabilizer group $\Gamma$,
then $Y:=E\times_F\spec\Omega$ is a projective scheme with very ample
coherent sheaf $\cE(k)\otimes\cO_Y$ that is generated by global sections,
and since $\Gamma$ acts faithfully on the fibers of $\cE(k)$
it follows that $\Gamma$ acts faithfully on $H^0(Y,\cE(k)\otimes\cO_Y)$.
Lemma \ref{l.flat} is proved.

\smallskip

Now Theorem \ref{thm.coq} is proved as follows.
Let $f\colon X \rightarrow F$ be a finite cover
of $F$ by a scheme (or algebraic space). By assumption
there is a surjection of sheaves ${\mathcal E} \rightarrow f_*{\mathcal
O}_X$. Let $V$ be the vector bundle associated with ${\mathcal E}$.
Then there is a closed immersion of $X$ into the stack $V$.
Since $X$ is representable and $X\rightarrow F$ is surjective,
$F$ is a quotient stack by Lemma \ref{l.vbqs}.

\begin{remark}
If, in the situation of Theorem \ref{thm.coq}, the stack $F$
admits a finite map to a scheme $Q$ (this occurs exactly when
$F$ has finite stabilizer and hence has a moduli space \cite{KM},
and the moduli space is a scheme) then
$Z$ (the algebraic space for which we have $F\simeq [Z/\GL_n]$)
is in fact a scheme.
The reason that $Z$ is a scheme is as follows:
Let $Y \rightarrow F$ be a finite cover of $F$ by a scheme.
Then, since $Z \rightarrow F$ is affine, the fiber product
$Z \times_F Y$ is an affine $Q$-scheme. Thus, by Chevalley's
theorem for algebraic spaces \cite[III.4.1]{Kn} it follows
that $Z$ is an affine $Q$-scheme as well.
\end{remark}

Finally, Theorem \ref{thm.dm} is a direct consequence
of Lemma \ref{l.vbqs},
provided we know that the tangent bundles
and higher jet bundles of smooth Deligne-Mumford stacks enjoy
faithful actions by the stabilizers of geometric points.

\begin{prop} Let $s,t\colon R\toto X$
be an \'etale presentation
of a smooth Deligne-Mumford stack $F$.
Let $\varphi\colon S \rightarrow X$ be the stabilizer group scheme. Assume
that no component of $S \smallsetminus e(X)$ dominates a component of $X$.
Then for some 
$k > 0$,
$S$ acts faithfully on the bundle
of $k$-jets in $X$.
\end{prop}

\begin{proof}
Let $x$ be a point in $X$.
Replacing $X$ by an \'etale cover if necessary, we may assume
the points of $\varphi^{-1}(x)$ all have residue field equal to the
residue field of $x$.
Then, for any $r\in \varphi^{-1}(x)$, $r\ne e(x)$,
the induced maps 
$$s^\#,t^\#\colon \widehat{{\mathcal O}}_{X,x}\rightarrow 
\widehat{{\mathcal O}}_{R,r}$$ are isomorphisms. 
Thus the composite
$$\widehat{{\mathcal O}}_{X,x} \stackrel{s^\#} \rightarrow
\widehat{{\mathcal O}}_{R,r} \stackrel{(t^{\#})^{-1}} \rightarrow 
\widehat{{\mathcal O}}_{X,x}$$ gives an automorphism
of the completed local ring $\widehat{{\mathcal O}}_{X,x}$.
By assumption on $S$, $s \neq t$ in a neighborhood
of $r \in R$ so the automorphism is nontrivial.
Thus, $r$ must act nontrivially on the vector space
${\mathcal O}_{X,x}/ {\mathfrak m}_x^k$
for some $k>0$.
Then, there exists $k$ such that the stabilizer group $\varphi^{-1}(x)$ acts
faithfully on
the space of $k$-jets at $x$.

By Noetherian induction on $X$, there is a $k$ for which the stabilizer
action on $k$-jets is faithful at all points of $X$.
\end{proof}

\begin{exml}
Let $k$ be a field of characteristic $p>0$.
The map $z\mapsto z^p-z$ realizes $C=\bP^1$ as a cyclic cover
of $\bP^1$, of degree $p$, branched only over infinity.
So, $\bP^1$ is the
coarse moduli space of the stack $F=[C/(\Z/p)]$, where a generator
of $\Z/p$ acts on $C$ by $[z:w]\mapsto [z+w:w]$.
The stabilizer of $F$ acts faithfully on the tangent bundle everywhere
except at the point over infinity.
For $n\ge 2$, the action on $n$-jets is faithful at all points.
\end{exml}
 
Finally, Corollary \ref{c.trivcenter} follows from the following
construction.
Let $F$ be a smooth Deligne-Mumford stack with finite stabilizer
$I_F\to F$.
There is an open dense substack $F^0$ on which the restriction
$I^0\to F^0$ of the stabilizer map is \'etale.
Let $J^0$ be the closure of $I^0$ in $I_F$; then
$\pi^0\colon J^0\to F$ is \'etale,
since any finite unramified morphism from a scheme to a normal
Noetherian scheme, such that every component dominates the target,
is \'etale.
Then $\cE:=\pi^0_*\cO_{J^0}$ is a
locally free coherent sheaf $F$.
We claim that the total space of the associated vector bundle has
trivial generic stabilizers, from which it follows (since $F$ embeds
in any vector bundle as the zero section) that $F$ is
a quotient stack.

Let $p\colon \spec\Omega\to F$ be a general geometric point of $F$,
with automorphism group $G$.
Then the fiber of $\pi^0$ over $p$ is canonically isomorphic to $G$,
and the action of $G$ over this fiber is by conjugation.
Since $G$, by hypothesis, has trivial center, the generic action on fibers of
$\pi^0_*\cO_{J^0}$ is faithful.

\subsection{A nonquotient stack}  We work out
Example \ref{e.noquotient}.
By Lemma \ref{l.vbqs}, if we can show the stack $F$
of Example \ref{e.noquotient} has no nontrivial vector bundles,
it follows that $F$ is not a quotient stack.

Let $Y=\spec \C[x,y,z]/(xy-z^2)$, with nonsingular locus
$Y^{\rm reg}$.
The nontrivial involution $i$ of $Y^{\rm reg}\times B(\Z/2)$ is
specified by
(it suffices to say how $i$ acts on pairs consisting of map $T\to Y^{\rm reg}$
and trivial $\Z/2$-torsor on $T$)
$$i(T\stackrel{f}\rightarrow Y^{\rm reg},\, T\times\Z/2\rightarrow T)=
(T\stackrel{f}\rightarrow Y^{\rm reg},\, f^*(
\A^2\smallsetminus\{0\}\rightarrow Y^{\rm reg})).$$
The stack $F$ is the union of two copies of $Y \times B(\Z/2)$,
glued via $i$.

By \cite{M-P} the scheme $Y$ has no nontrivial vector bundles
(this fact holds more generally for any affine toric variety
\cite{Gu}),
and since $Y$ is normal and the glueing is over a locus whose
complement has codimension 2, the scheme $Y\amalg_{Y^{\rm reg}}Y$
(this is the scheme over which $F$ is a gerbe)
also has no nontrivial vector bundles.
Every vector bundle on $F$ splits into $(+1)$- and $(-1)$-eigenbundles
for the stabilizer action,
so we deduce that the $(+1)$-eigenbundle is trivial.

We claim the $(-1)$-eigenbundle is zero.
Let ${\mathcal F}$ be a locally free coherent sheaf on $F$ such that the
stabilizer action is multiplication by $-1$ on sections,
and let ${\mathcal F}_i$ ($i=1, 2$) denote the restriction of ${\mathcal F}$
over the $i$-th copy of $Y$.
Then there is a given isomorphism
$${\mathcal F}_1|_{F^{\rm reg}}\simeq i^*({\mathcal F}_2|_{F^{\rm reg}}).$$
Let ${\mathcal L}$ denote the (pullback to $F^{\rm reg}$ of the) 
unique 2-torsion invertible sheaf on $Y^{\rm reg}$;
for any locally free coherent sheaf ${\mathcal F}$ on $F^{\rm reg}$ such that
the stabilizer acts by $(-1)$ we have
$i^*{\mathcal F}\simeq
{\mathcal L}\otimes {\mathcal F}$.
Hence $${\mathcal F}_1|_{Y^{\rm reg}}\simeq
{\mathcal L}\otimes ({\mathcal F}_2|_{Y^{\rm reg}}).$$
But this is impossible unless ${\mathcal F}_1={\mathcal F}_2=0$,
for otherwise ${\mathcal F}_i|_Y$ ($i=1,2$) is free of some rank $m\ge 1$,
and hence we have ${\mathcal L}^{\oplus m}\simeq {\mathcal O}^{\oplus m}$
on $Y^{\rm reg}$.  But $Y^{\rm reg}$ sits inside $[\A^2/(\Z/2)]$ with
complement of codimension 2,
so this implies an isomorphism
on $[\A^2/(\Z/2)]$ between a
free coherent sheaf and a nontrivial locally free coherent sheaf.

\subsection{Gerbes and the Brauer group}
Here we prove Theorem \ref{t.brauer}.
For (i) $\Rightarrow$ (iii),
let $\gamma\in H^1(X,PGL_n)$ be the class of an Azumaya algebra
representing a given cohomological Brauer group element
$\beta\in H^2(X,\G_m)$.
If $P\to X$ is the $PGL_n$-bundle associated with $\gamma$,
then by the definition of the boundary map in nonabelian cohomology,
the gerbe represented by $\beta$ is $[P/GL_n]$.

For (iii) $\Rightarrow$ (i), we note that if $F$ is the $\G_m$-gerbe
associated with $\beta$, then a vector bundle $B$ on $F$ decomposes into
eigenbundles indexed by characters in $\widehat \G_m$.
Given a faithful stabilizer action on fibers, the characters whose
eigenbundles are nonzero must generate $\widehat \G_m$;
then the decomposition of $B^{\otimes r}\otimes (B^\vee)^{\otimes s}$
for suitable integers $r$ and $s$ has nonzero eigenbundle $B_1\to F$
for the unit character.
The complement of the zero section of $B_1$ is a
Brauer-Severi scheme over $X$,
and the associated Azumaya algebra represents $\beta$.

As we have remarked, (i) $\Rightarrow$ (ii) is well-known.
The implication (ii) $\Rightarrow$ (iii) is an immediate consequence
of Lemma \ref{l.flat}.

Finally, suppose $\beta$ is $n$-torsion with pre-image
$\alpha$ in $H^2(X,\mu_n)$, and let us show
(iii) $\Leftrightarrow$ (iv).
Let $F$ be the $\G_m$-gerbe associated with $\beta$,
and let $F'$ be the $\mu_n$-gerbe associated with $\alpha$.
There is a natural representable morphism $\pi\colon F'\to F$,
hence (iii) $\Rightarrow$ (iv).
For the reverse implication, let $\cE$ be a locally free coherent
sheaf on $F'$ such that the stabilizer action on sections is faithful.
Consider the quasi-coherent sheaf $\cF:=\pi_*\cE$ with its decomposition
into eigensheaves $\cF=\bigoplus_{\lambda\in\widehat\G_m}\cF_\lambda$.
We claim each $\cF_\lambda$ is locally free of finite type, and
$\cF_\lambda\ne 0$
if and only if the eigensheaf of $\cE$ corresponding to
the restriction of $\lambda$ to $\mu_n$ is nonzero.
Indeed, it suffices to verify the claims \'etale locally,
and the claims hold in the case of trivial gerbes.
We can choose a finite set $S$ of characters which generates $\G_m$ such that
$\cF_\lambda\ne 0$ for every $\lambda\in S$.
Then $\bigoplus_{\lambda\in S}\cF_{\lambda}$ is
a locally free coherent sheaf on $F$,
such that the stabilizer action on sections is faithful.

\end{document}